\documentclass{amsart}
\usepackage{amssymb, latexsym,graphicx,bm,amscd,amsthm, amsmath, latexsym, cite, subfigure, multirow}
\usepackage[all]{xy}
\newcommand{\integer}[1]{\left[ #1 \right]}

\newtheorem{theorem}{Theorem}
\newtheorem{proposition}{Proposition}[section]
\newtheorem{lemma}{Lemma}[proposition]
\newtheorem{corollary}{Corollary}[proposition]
\newtheorem*{theorem-mug}{Theorem 1, \cite{ralston2}}
\theoremstyle{remark}
\newtheorem{example}{Example}[proposition]

\title{$1/2$-Heavy sequences driven by rotation}
\author{David Ralston}
\address{Department of Mathematics, Ben Gurion University of the Negev \\POB 653
 Beer Sheva, 84105 \\ Israel.}
\date{\today}
\email{ralston.david.s@gmail.com}
\begin{document}
\begin{abstract}
We investigate the set of $x \in S^1$ such that for every positive integer $N$, the first $N$ points in the orbit of $x$ under rotation by irrational $\theta$ contain at least as many values in the interval $[0,1/2]$ as in the complement.  By using a renormalization procedure, we show both that the Hausdorff dimension of this set is the same constant (strictly between zero and one) for almost-every $\theta$, and that for every $d \in [0,1]$ there is a dense set of $\theta$ for which the Hausdorff dimension of this set is $d$.
\end{abstract}
\maketitle
\section{Introduction}
Let $X = \mathbb{R}/\mathbb{Z}$ be the unit circle, with addition modulo one, and for irrational $\theta$ fixed define
\[f(x)=\chi_{[0,1/2]}(x)-\chi_{(1/2,1)}(x), \quad T(x) = x + \theta \bmod 1, \quad S_n(x) = \sum_{i=0}^{n-1}f \circ T^i (x).\]
We will explicitly present the sets
\[H_{\theta}=\left\{ x \in [0,1) : S_n(x) \geq 0, \quad n=1,2,\ldots\right\},\]
\[H^*_{\theta} =\left\{ x \in [0,1) : S_n(x) > 0, \quad n=1,2,\ldots\right\}.\]
These sets are known as the \textit{heavy} and \textit{strictly heavy} sets, respectively.  The sequence $S_n(x)$ is the $1/2$-\textit{discrepancy sequence} defined by a given $x$ and the rotation parameter $\theta$.  It follows from \cite{MR931186} that $H_{\theta}$ is nonempty and from \cite{springerlink:10.1007/BF01896805} that it is of measure zero.  The special case of $x=\theta$ was studied in \cite{MR2515388} (where $S_n(\theta)\geq 0$ was termed `\textit{$\{n \theta\}$ is a $1/2$-heavy sequence}'), and Theorem 1 in that work may be interpreted as saying that $0 \in H^*_{\theta}$ if and only if $\theta$ has a particular expansion as a continued fraction: all partial quotients of odd index are themselves even.

\begin{theorem}\label{theorem - strictly heavy singleton}
For every irrational $\theta$, $H^*_{\theta}$ is a singleton.
\end{theorem}
\begin{corollary}
$H^*_{\theta}=\{0\}$ if and only if all partial quotients of $\theta$ of odd index are themselves even.
\end{corollary}
\begin{theorem}\label{theorem - heavy constant}
There is some constant $c \in (0,1)$ such that for almost-every $\theta$, $\dim_H\left( H_{\theta} \right) = c$.
\end{theorem}

\begin{theorem}\label{theorem - heavy dense}
Given any $d \in [0,1]$, there is a dense set of $\theta$ for which $\dim_H \left( H_{\theta} \right) =d$.
\end{theorem}

Several examples making explicit use of the techniques developed are presented.

\section{The Renormalization Procedure}\label{section - renormalization}
We will freely use standard continued fraction notation.  For $\theta \in (0,1)\setminus \mathbb{Q}$, let
\[ \theta=[a_1,a_2,a_3, \ldots]= \cfrac{1}{a_1+\cfrac{1}{a_2+\cfrac{1}{\ddots}}}.\]
The first of the following inequalities is standard, and the second is easy to verify:
\begin{subequations}
\begin{align}
\frac{1}{(2n+1)(m+1)} < &1-2n\theta < \frac{1}{2nm} & (a_1=2n, \, a_2=m) \label{eqn - estimate even}\\
\frac{1}{2n+1} < &1-2n\theta < \frac{1}{n+1}& (a_1=2n+1 \neq 1) \label{eqn - estimate odd}
\end{align}
\end{subequations}
The \textit{Gauss map} is given by
\[\gamma(\theta) = \frac{1}{\theta} \bmod 1, \quad \gamma([a_1,a_2,\ldots])=[a_2,a_3,\ldots].\]
Define the function $g$ by
\begin{equation}
g(\theta) = \begin{cases} 1-\theta=[a_2+1,a_3,\ldots] & (a_1(\theta)=1)\\ \frac{1}{1+\gamma(\theta)}=[1,a_2,a_3,\ldots]& (a_1(\theta) = 2n+1 \neq 1) \\ \gamma^2(\theta)=[a_3,\ldots] & (a_1(\theta)=2n)\end{cases}
\label{eqn - g}
\end{equation}

As rotation by any irrational $\theta$ is topologically minimal, for any proper interval $I' \subset [0,1)$, we define for every $x$
\[n(x) = \min\left\{ n \in \mathbb{N} : T^n(x) \in I' \right\},\]
and the \textit{induced map} on $I'$ is defined for $x \in I'$ by $T_{I'}(x) = T^{n(x)}(x)$.  Let $\integer{x} = x - \left( x \bmod 1 \right)$ denote the \textit{integer part} or \textit{floor} of $x$.  The following information may all be derived from \cite[\S 3]{ralston}, where the specified renormalization is studied in detail:
\begin{proposition}\label{proposition - induction scheme}
For a fixed $\theta \in (0,1) \setminus \mathbb{Q}$, let $I'$ be the interval $[0,\delta)$, where
\[\delta = 1 - 2 \integer{\frac{a_1}{2}}\theta.\]  Then if the endpoints of $I'$ are identified with one another, the first return map $T_{I'}$ is rotation by $g(\theta) \delta$.

For those $\theta<1/2$ (so that $I' \neq [0,1)$), for all $x \in [0,\delta/2]$ we have 
\[1 \leq i \leq n(x) \quad \Longrightarrow \quad S_i(x) \geq 1,\] and $S_{n(x)}(x) =1$.

For those $x \in (\delta/2,\delta)$, on the other hand,
\[1 \leq i \leq n(x) \quad \Longrightarrow \quad S_i(x) \geq -1,\] and $S_{n(x)}(x) = -1$.
\end{proposition}

This ergodic properties of this renormalization map have also been presented:
\begin{theorem-mug}
There exists a unique probability measure $\mu_g$ on the circle, mutually absolutely continuous with respect to Lebesgue measure, such that $\mu_g$ is $g$-invariant and ergodic.  The system $\{S^1,\mu_g,g\}$ is in fact exponentially CF-mixing, and both Radon-Nikodym derivatives $d\mu_g/dx$ and $dx/d\mu_g$ are essentially bounded.
\end{theorem-mug}

\section{Structure of Heavy Sets}\label{section - heaviness}

Denote the Hausdorff, upper Minkowski box, and lower Minkowski box dimensions of a set $S$ by $\dim_{H}(S)$, $\dim_{B}(S)$, and $\dim_b(S)$, respectively, so for any set $S$ we have \[\dim_H(S) \leq \dim_b(S) \leq \dim_B(S).\]

Let $E_0$ be a closed, connected interval in $[0,1)$, and define a sequence of sets
\[E_0 \supset E_1 \supset E_2 \supset \ldots\]
by requiring that each $E_i$ be a disjoint union of closed intervals and that the maximal length of an interval in $E_i$ tends to zero as $i \rightarrow \infty$.  Define $F = \cap_{i=0}^{\infty} E_i$.  Suppose further that each $E_i$ contains at least $m_{i+1}$ intervals of $E_{i+1}$, and that the intervals of $E_i$ are separated by at least $\epsilon_i>0$, where $\epsilon_{i+1}<\epsilon_i$.  Then by \cite[Example 4.6]{MR2118797}:
\begin{equation}\label{equation - hdim falconer example}
\dim_H(F) \geq \liminf_{k \rightarrow \infty} \frac{\sum_{i=1}^{k} \log m_i}{-\log m_{k+1} - \log \epsilon_{k+1}}.
\end{equation}

A sequence of numbers $\{m_i\}$ will be said to be \textit{not exceptionally irregular} if for every $\epsilon>0$, there is some number $N$ such that for all $k \geq N$, we have
\begin{equation}\label{eqn - not exceptionally irregular}
m_k < \prod_{i=1}^{k-1} m_i ^{\epsilon}.
\end{equation}

\begin{proposition}\label{proposition - hausdorff is box}
Suppose that each interval of $E_i$ contains exactly $m_{i+1} \geq 1$ intervals of $E_{i+1}$, each of which is of length exactly $\delta_{i}\leq 1$ times the length of the intervals in $E_i$, and furthermore suppose that the intervals of $E_i$ are separated by gaps at least as large as the intervals themselves, and this length converges to zero as $i \rightarrow \infty$.  Finally, assume that the $m_k$ are not exceptionally irregular \eqref{eqn - not exceptionally irregular}.  Then
\[\dim_H(F) = \dim_b(F) = \liminf_{k \rightarrow \infty} \frac{\sum_{i=1}^k \log m_i}{-\sum_{i=1}^{k} \log \delta_i}.\]
\begin{proof}
First, note that all intervals in $E_i$ are of length $\delta_1 \cdots \delta_i$ times the length of $E_0$, and this length also serves as $\epsilon_i$ in \eqref{equation - hdim falconer example}.  From \eqref{equation - hdim falconer example}, the definition of lower Minkowski box dimension $\dim_b$ and the fact that $\dim_H \leq \dim_b$, we have that
\[\liminf_{k \rightarrow \infty}\frac{\sum_{i=1}^k \log m_i}{-\log m_{k+1} - \sum_{i=1}^k \log \delta_i} \leq \dim_H(F) \leq \dim_b(F) = \liminf_{k \rightarrow \infty}\frac{\sum_{i=1}^k \log m_i}{-\sum_{i=1}^k \log \delta_i},\]
and the condition \eqref{eqn - not exceptionally irregular} exactly ensures that equality is forced.
\end{proof}
\end{proposition}

Returning our attention to the sets $H_{\theta}$, let $\theta$ be fixed.  Note that we trivially have $H^*_{\theta} \subset H_{\theta} \subset [0,1/2]$.
\begin{lemma}\label{lemma - reversal of heavy}
If $a_1(\theta)=1$, then (with subtraction defined pointwise)
\[H_{\theta} = \frac{1}{2} - H_{g(\theta)},\quad H^*_{\theta} = \frac{1}{2} - H^*_{g(\theta)}.\]
\begin{proof}
$g(\theta)=1-\theta$ for $a_1(\theta)=1$, and for any $\theta$, the reader may verify that
\[x+i \theta \leq \frac{1}{2} \quad \Longleftrightarrow \quad (1/2-x)+i(1-\theta) \leq 1/2.\qedhere \]
\end{proof}
\end{lemma}

\begin{lemma}\label{lemma - rescaled heavy set}
Let $\theta$ be fixed with $a_1(\theta) \neq 1$, and $\delta=1-2[a_1/2]\theta$.  Then
\[H_{\theta} \cap [0,\delta) = \delta \cdot H_{g(\theta)},\] where multiplication is defined pointwise.
\begin{proof}
Suppose $x \in H_{\theta}\cap[0,\delta)$.  Then as all ergodic sums are nonnegative, certainly the ergodic sums along the subsequence of times $n_i(x)$, the $i$-th return to $[0,\delta)$, are also nonneagtive.  However, by Proposition \ref{proposition - induction scheme} we have
\[S_{n_i}(x) = \sum_{j=0}^{i-1} \left(\chi_{[0,\delta/2]}-\chi_{(\delta/2,\delta)}\right)(T_{I'}^j (x)),\] as the cumulative sums through return to $I'$ are exactly $\pm 1$ on these intervals.  As this function is simply $f(\delta x)$, and $T_{I'}$ is rotation by $g(\theta)\cdot \delta$ (recall that the endpoints of $I'$ are identified with one another), we have 
\[H_{\theta}\cap[0,\delta] \subset \delta \cdot H_{g(\theta)}.\]
For the other containment, we note (again using Propositon \ref{proposition - induction scheme}) that the cumulative sum $S_{n(x)}(x)$ is in fact the minimal value reached until return to $n(x)$, so if
\[\sum_{j=0}^{k-1} \left(\chi_{[0,\delta/2]}-\chi_{(\delta/2,\delta)}\right)(T_{I'}^j (x)) \geq 0\] for $k=0,1,\ldots,i$, then $S_n(x) \geq 0 $ for all $n=0,1,\ldots,n_i(x)$.
\end{proof}
\end{lemma}
\begin{corollary}\label{corollary - empty bit}$H_{\theta} \cap (\delta/2,\delta)=\emptyset$.
\begin{proof}$H_{g(\theta)}\subset [0,1/2]$.
\end{proof}
\end{corollary}
\begin{lemma}\label{lemma - critical cutoff}
Let $a_1(\theta) \neq 1$.  Then $H_{\theta} \cap (\delta/2+\theta,1) = \emptyset$.
\begin{proof}
Clearly no $x>1/2$ is in $H_{\theta}$.  So assume that $a_1(\theta)=2n$ or $2n+1$, but in either case $n \neq 0$.  If $n=1$, then $\delta/2+\theta = 1/2$, so without loss of generality assume $n \geq 2$.  Let $x \in (1/2-(n-1)\theta,1/2]$.  As $x>1/2-(n-1)\theta$, evaluating $f$ along the orbit of $x$ produces no more than $n-1$ consecutive values of $+1$.  Since $n\theta<1/2$, however, after this initial string of positive values, the orbit includes at least $n$ consecutive values of $-1$, producing a negative sum in the orbit of $x$.
\end{proof}
\end{lemma}

\begin{lemma}\label{lemma - ignorable region}
Suppose $a_1(\theta)\neq 1$.  Then for $x \in [\delta, \delta/2+\theta]$, 
\[x \in H_{\theta} \quad \Longleftrightarrow \quad (x+2n \theta) \in H_{\theta}.\]
\begin{proof}
Let $a_1(\theta)=2n$ or $2n+1$; in either case $n \neq 0$.  As \[\delta/2=1/2-n \theta < x < 1/2-(n-1)\theta,\] the orbit begins with exactly $n$ consecutive values within $[0,1/2]$.  As in Lemma \ref{lemma - critical cutoff}, it then contains at least $n$ consecutive values within $(1/2,1)$.  So the first $2n$ values of $f$ along the orbit of $x$ consist of exactly $n$ consecutive values of $+1$ followed by $n$ consecutive values of $-1$.  The ergodic sums have not been negative before time $2n$, and for $k \geq 0$ we have
\[S_{2n+k}(x) = S_{2n}(x)+S_k(x+2n \theta)=S_k(x+2n\theta). \qedhere\]
\end{proof}
\end{lemma}
\begin{corollary}\label{corollary - strictly heavy cutoff}
$H^*_{\theta} \subset [0,\delta/2]$.
\begin{proof}
We have shown that for $x>\delta/2$, there is some $n$ for which $S_n(x) \leq 0$.
\end{proof}
\end{corollary}
\begin{lemma}\label{lemma - counting islands}
If $a_1(\theta)=2n+1 \neq 1$, then
\[\frac{\delta}{2}+\delta > \frac{1}{2}-(n-1)\theta.\]
If $a_1(\theta)=2n$ and $a_2(\theta)=a_2$, then
\[\frac{\delta}{2}+a_2 \delta < \frac{1}{2}-(n-1)\theta< \frac{\delta}{2} + (a_2+1)\delta.\]
Furthermore, if $a_1(\theta)=2n$, then
\[(a_2+1)\delta < \frac{1}{2}-(n-1)\theta\] if and only if $a_3(\theta)=1$.
\begin{proof}
All of these facts are simple computatations using elementary inequalities of continued fractions.  We include only the computations behind the last statement to give an indication of the techniques necessary, where each line is equivalent to the next:
\begin{align*}
(a_2+1)\delta &< \frac{1}{2}-(n-1)\theta \\
a_2 \left( \frac{1}{\theta}-2n\right)+ \frac{1}{\theta}-2n &< \frac{1}{2}\left( \frac{1}{\theta}-2n\right) +1\\
\gamma(\theta) \left(a_2 +1\right) &< \frac{\gamma(\theta)}{2}+1\\
\frac{1}{2} &< \frac{1}{\gamma(\theta)}-a_2\\
\frac{1}{2} &< \gamma^2(\theta).\qedhere
\end{align*}
\end{proof}
\end{lemma}

\begin{figure}\centering{\fbox{
\xymatrix@!0@C=1.6pc@R=2pc{
&&&&&&&&&&&& &\frac{\delta}{2}+a_2\delta \ar[d]\\
&&&&&&&&&&&&&{\bullet}\ar@{.}[r]&{\bullet} \\
{\bullet}\ar@{-}[rr] &&{\bullet}\ar@{.}[rr]&& {\bullet}\ar@{-}[rr]& \ar@/_1pc/[llll]&{\bullet}\ar@{.}[rr]&\ar@/_1pc/@{.>}[llll]&{\bullet}\ar@{-}[rr]&\ar@/_1pc/[llll]_{\textit{etc.}}&{\bullet}\ar@{.}[r] & \ar[urr] \ar[drr] &&&\frac{\delta}{2}+\theta \ar[u]_{a_3\neq 1} \ar[d]^{a_3=1} \\
0&&\frac{\delta}{2}&&\delta &&\delta + \frac{\delta}{2} && 2 \delta && \frac{\delta}{2}+2 \delta &&& {\bullet}\ar@{-}[r]&{\bullet}\\
&&&&&&&&&&&&&(a_2+1) \delta \ar[u] 
}}}
\caption{The `proof by picture' of Corollary \ref{corollary - empty bit} in the case $a_1=2n$.  The solid intervals intersect $H_{\theta}$ with affine images of $H_{g \theta}$ by iterated application of Lemma \ref{lemma - ignorable region} (indicated by the left-pointing arrows) combined with Lemma \ref{lemma - rescaled heavy set}, while the dotted intervals do not intersect $H_{\theta}$ by iterating Lemma \ref{lemma - ignorable region} and considering Corollary \ref{corollary - empty bit}.  In the event that $a_3=1$, there is an additional subset of an affine image of $H_{g \theta}$ in the interval $[(a_2+1)\delta,1/2-(n-1)\theta]$.  Placement of the cutoff point $\delta/2+\theta$ is determined by Lemma \ref{lemma - counting islands}.}\label{figure - proof by picture}
\end{figure}
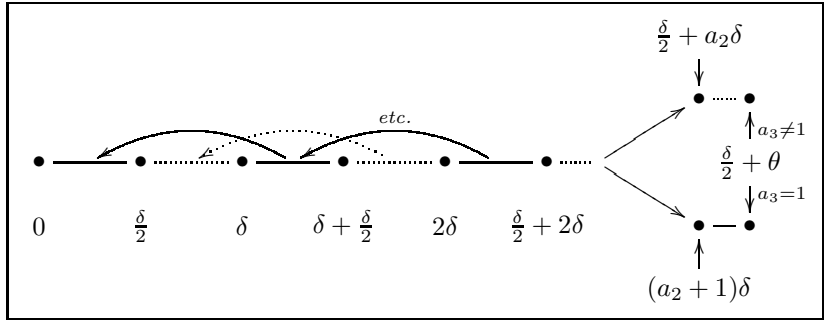

\begin{corollary}\label{corollary - heavy decomposition}
For all $\theta$ such that $a_1(\theta) \neq 1$, $H^*_{\theta}=\delta \cdot H^*_{g(\theta)}$.
If $a_1(\theta)=2n+1 >1$, then (with addition and multiplication defined pointwise),
\begin{itemize}
\item $H_{\theta} \subset [0,\delta/2] \cup [\delta,1/2-(n-1)\theta]$,
\item $H_{\theta} \cap [0,\delta/2] = \delta \cdot H_{g(\theta)}$,
\item $H_{\theta} \cap [\delta,1/2-(n-1)\theta] = \left(H_{\theta}\cap[0,(n+1)\theta-1/2]\right)+\delta$.
\end{itemize}
If $a_1(\theta)=2n$ and $a_3(\theta)\neq 1$,
\begin{itemize}
\item $H_{\theta} \subset \bigcup_{i=0}^{a_2} \left([0,\delta/2]+i\delta\right)$, for $i=0,1,\ldots,a_2$,
\item $H_{\theta} \cap \left( [0,\delta/2]+i\delta\right) = \delta \cdot H_{g(\theta)}+i \delta$,
\end{itemize}
If $a_1(\theta)=2n$ and $a_3(\theta)=1$,
\begin{itemize}
\item  $H_{\theta} \subset [(a_2+1)\delta,1/2-(n-1)\theta] \cup \left( \bigcup_{i=0}^{a_2} \left([0,\delta/2]+i\delta\right)\right)$,
\item $H_{\theta} \cap \left( [0,\delta/2]+i\delta\right) = \delta \cdot H_{g(\theta)}+i \delta$, for $i=0,1,\ldots,a_2$
\item $H_{\theta} \cap [(a_2+1)\delta, 1/2-(n-1)\theta] = \left(H_{\theta}\cap [0,(n+1)\theta-1/2]\right)+(a_2+1)\delta$,
\end{itemize}
\begin{proof}
The proof is a combination of all previous discussion and is left to the reader.  See Figure \ref{figure - proof by picture} for a `proof by picture' which illustrates all the necessary steps for the case $a_1=2n$.
\end{proof}
\end{corollary}

\section{Proof of Theorem \ref{theorem - strictly heavy singleton}}
Let $\theta_i=g^i(\theta)$, $\delta_i = 1-2[a_1(\theta_i)/2]\theta_i$, and
\[p_i = \#\left\{ j = 0,1,\ldots, i-1: a_1(\theta_i)=1\right\} \bmod 2,\]
where $\#$ denotes cardinality and $p_0=0$.

We trivially know that $H^*_{\theta} \subset [0,1/2]$, so denote $E^*_0=[0,1/2]$ and $I_0 = [0,1)$.  Now, inductively define $I_{i+1}$ and $E^*_{i+1}$, using $E^*_i=[a,b]$.  If $p_{i+1}=0$, let \[I_{i+1}=[a,a + \delta_i |I_i| ), \, E^*_{i+1}=[a,a+\delta_i |I_i|/2],\] where $|I_i|$ is the length of $I_i$, which is easily inductively seen to be $2(b-a)$.  If $p_{i+1}=1$, let 
\[I_{i+1}=(b-\delta_i|I_i|,b], \, E^*_{i+1} =[b-\delta_i |I_i|/2,b].\]  It follows from Proposition \ref{proposition - induction scheme} (reversals of orientation induced by $a_1(\theta_i)=1$ are accounted for by $p_i$) that ergodic sums with respect to the function $\chi_{E^*_{i}}(x)-\chi_{I_i \setminus E^*_{i}}(x)$ and the transformation $T_{I_i}$ are never less than one through the return time to $I_{i+1}$ exactly on the set $E^*_{i+1}$, which is the total ergodic sum at this time, and on the interval $I_{i+1}\setminus E^*_{i+1}$ the ergodic sums are never less than minus one through their return, and this is the value of the ergodic sums at the return time.  Furthermore, on $I_i \setminus E^*_{i+1}$ there is some non-positive ergodic sum (Corollary \ref{corollary - strictly heavy cutoff}).  Inductively, then, $E^*_{i}$ is exactly the set of points in $I_0$ for which $S_n(x) \geq 1$ for all $n \leq n(x)$, the return time to $I_i$.

The length of $I_i$ is given by $\delta_0 \delta_1 \cdots \delta_{i}$.  For $a_1(\theta_i)\neq 1$ (i.e. $\theta_i < 1/2$), we have $\delta_i < 1/2$, and since $g:(1/2,1) \rightarrow (0,1/2)$, of any consecutive $\delta_i$, $\delta_{i+1}$, at least one is less than one half.  This shows both that the intersection of all $E^*_i$ is a singleton and the set $H^*_{\theta}$ (the return time to $I_i$ must diverge as the length decreases to zero).

\begin{corollary}\label{corollary - always infinite}
The set $H_{\theta}$ is always infinite.
\begin{proof}
We must have infinitely many $N_i$ such that $S_{N_i}(H^*_{\theta})=1$.  Otherwise, the largest such $N_i$ would yield $H^*_{\theta}+N_i\theta$ as another strictly heavy point (and the set of such $N_i$ is not empty, containing $1$).  For each such $N_i$ we see 
\[S_{k}(H^*_{\theta}+N_i\theta) = S_{N_i+k}(H^*_{\theta})-S_{N_i}(H^*_{\theta}) \geq 1-1=0. \qedhere\]
\end{proof}
\end{corollary}

That $S_n(x)$ form an additive cocycle allows us to immediately conclude (for $a_1(\theta)\neq 1$, i.e. $\theta<1/2$):
\begin{equation}\label{eqn - partial copy singletons}
H_{\theta} \cap [\theta,1/2] = H^*_{\theta}+\theta.
\end{equation}

\section{Proof of Theorem \ref{theorem - heavy dense}}
Before we construct the sets $E_i$ to determine $\dim_H(H_{\theta})$, we refine Corollary \ref{corollary - heavy decomposition} using Theorem \ref{theorem - strictly heavy singleton}:

\begin{lemma}\label{lemma - ignore partial copies}
If $a_1(\theta)=2n+1\neq 1$, then 
\[H_{\theta} \cap [\delta, 1/2]\]
is a singleton (and therefore an isolated point within $H_{\theta}$).  If $a_1(\theta)=2n$ and $a_3(\theta)=1$, then 
\[H_{\theta} \cap [(a_2(\theta)+1)\delta, 1/2-(n-1)\theta]\] is a singleton (and therefore an isolated point within $H_{\theta}$).
\begin{proof}
First let $a_1(\theta)=2n+1$.  In this case $\delta/2<\theta<\delta$, and  $H_{\theta} \cap [\theta,1/2]$ is a singleton  by \eqref{eqn - partial copy singletons}.  We already know that $H_{\theta} \cap (\delta/2,\delta) = \emptyset,$ so that singleton is in $[\delta,1/2]$.
In the event that $a_1(\theta)=2n$ and $a_1(\theta)=1$, we verify
\[ \renewcommand\arraystretch{1.5}\begin{array}{r c c c l}
\frac{1}{2} &< &\gamma^2(\theta) &< &1\\
\frac{1}{2}+a_2(\theta) &<& \frac{1}{\gamma(\theta)} &<& 1 + a_2(\theta)\\
\left( \frac{1}{\theta} - 2n\right) \left( \frac{1}{2} + a_2 \right) &<& 1 &<& \left(\frac{1}{\theta}-2n \right)\left(1+a_2 \right)\\
\delta \left(\frac{1}{2}+a_2\right) &<& \theta &<& \delta \left(1+a_2\right).
\end{array}\]
So, $H_{\theta}\cap [\theta,1/2]$ is a singleton within $[(a_2+1)\delta,1/2-(n-1)\theta]$; from Corollary \ref{corollary - heavy decomposition} we have
\[H_{\theta} \cap \left( \delta \left( \frac{1}{2} + a_2 \right), \delta (1+a_2)\right) = \emptyset. \qedhere\]
\end{proof}
\end{lemma}

We now let $E_0=[0,1/2]$, and $\theta_i$, $\delta_i$, $p_i$ are defined as in the proof of Theorem \ref{theorem - strictly heavy singleton}.  We further define $S_0=\emptyset$, $\Delta_0=1$, and
\[\Delta_i = \prod_{j=0}^{i-1} \delta_i.\]
Assume that $E_i$ is given by a collection of closed intervals, each of which are of length $(1/2)\Delta_i$ and which are separated from one another by at least $(1/2)\Delta_i$.  If $a_1(\theta_i)=1$, then if we let $E_{i+1}=E_i$ and $S_{i+1}=S_i$, this assumption will be maintained, as for $\theta_i>1/2$ we have $\delta_i=1$.

Inductively assume that for each interval $[a,b] \in E_i$, we have
\[H_{\theta} \cap [a,b] = \Delta_i \cdot H_{\theta_i}+a,\] which is clearly true for $i=0$.  Then for each such $[a,b] \in E_i$ we add to $E_{i+1}$ a collection of intervals according to several cases, given in Table \ref{table - decomposition}.

\begin{table}[hbt]
\begin{tabular}{|c| c| c|}
\hline \multicolumn{2}{|c|}{ Case} & Interval(s) in $E_{i+1}$\\
\hline \multicolumn{2}{|c|}{$a_1(\theta_i)=1$} & $[a,b]$ \\
\hline
\multirow{2}{*}{$a_1(\theta_i)=2n+1\neq1$}
	&$p_i=0$ & $[a,a+(1/2)\Delta_{i+1}]$\\
	&$p_i=1$ & $[b-(1/2)\Delta_{i+1},b]$\\
\hline
\multirow{2}{*}{$a_1(\theta_i)=2n, \, a_2(\theta)=a_2$}
	&$p_i=0$ & $[a,a+(1/2)\Delta_{i+1}]+j\Delta_{i+1}$ \, $(j=0,\ldots,a_2)$\\
	&$p_i=1$ & $[b-(1/2)\Delta_{i+1},b] - j\Delta_{i+1}$\, $(j=0,\ldots,a_2)$\\
\hline
\end{tabular}
\caption{The decomposition of $E_i$ into $E_{i+1}$ according to $\theta_i$.}\label{table - decomposition}
\end{table}

So each interval in $E_{i+1}$ is of length $(1/2)\Delta_{i+1}$ and separated from other intervals in $E_{i+1}$ by the same amount, by applying Corollary \ref{corollary - heavy decomposition} and Proposition \ref{proposition - induction scheme}, we know that each such interval contains the desired affine image of $H_{\theta_{i+1}}$.  The terms $p_{i}$ track the parity of the number of changes of orientation caused by applying Lemma \ref{lemma - reversal of heavy} as in the proof of Theorem \ref{theorem - strictly heavy singleton}.

Finally, if $a_1(\theta_i)=2n+1$ or $a_1(\theta_i)=2n$, $a_3(\theta_i)=1$, we define $S_{i+1}$ to be the union of $S_i$ along with the (finitely many, possible zero) isolated points which arise in those smaller intervals listed in Corollary \ref{corollary - heavy decomposition} (applying Lemma \ref{lemma - ignore partial copies}).

\begin{proposition}
$H_{\theta} = E \cup S$, where
\[E = \bigcap_{i=0}^{\infty}E_i, \qquad S= \bigcup_{i=0}^{\infty} S_i.\]
\begin{proof}
First, let $x \in S_i$, where $i$ is the minimal index such that $x \in S_i$.  Then $x$ was added to $S_i$ precisely because it is the affine image of the point $H^*_{\theta_i}$ within $H_{\theta_{i}} \subset [a,b] \in E_{i-1}$, and our inductive hypothesis ensures that $x \in H_{\theta}$.  Similarly, if $x \in E$, then $x \in E_i$ for all $i$.  Just as in the proof of Theorem \ref{theorem - strictly heavy singleton}, then, $x \in H_{\theta}$; nonnegative sums are maintained through return times to an interval whose length is converging to zero.

Suppose now that $x \notin (E \cup S)$.  Our inductive construction of the $E_i$ and $S_i$ then guarantees that $x \notin H_{\theta}$: by applying Corollary \ref{corollary - heavy decomposition}, it must fail to be in any of the affine images of $H_{\theta_i}$ that comprise the set $H_{\theta}$.
\end{proof}
\end{proposition}

Note that $S$ is a countable set, as each $S_i$ is finite (alternately, from Lemma \ref{lemma - ignore partial copies} we know that all points in $S$ are separated from each other).  Therefore, $\dim_H(H_{\theta})=\dim_H(E)$.  Define the functions
\[f_1(\theta) = 1-2\integer{\frac{a_1(\theta)}{2}}\theta,\]
\[f_2(\theta) = \begin{cases} 1 & (a_1(\theta)=1 \bmod 2)\\ a_2(\theta)+1 & (a_1(\theta)=0 \bmod 2). \end{cases}\]
Then the set $E_{i+1}$ is formed from $E_i$ by letting each interval in $E_i$ be replaced with $f_2(\theta_i)$ intervals of length scaled by $f_1(\theta_i)$, all of which are separated from one another by gaps at least as large as the intervals themselves.

\begin{lemma}\label{lemma - no fast growth}
For almost every $\theta$, the sequence $\{f_2(\theta_i)\}$ is not exceptionally irregular \eqref{eqn - not exceptionally irregular}.
\begin{proof}
Recall that $\theta_i=g^i(\theta)$, and $g$ is ergodic with respect to the continuous measure $\mu_g$.  By the Birkhoff ergodic theorem, then, for almost every $\theta$ we have
\[\lim_{n \rightarrow \infty} \frac{1}{n} \sum_{i=0}^{n-1}\log (f_2(\theta_i)) = \int_{S^1} \log(f_2(\theta))d \mu_g,\]
provided that $\log(f_2) \in L^1(S^1,\mu_g)$.  Since $\mu_g$ is continuous, however, it suffices to show that $\log(f_2) \in L^1(S^1,dx)$.  Note that $f_2$ is constant on intervals of $\theta$ whose continued fraction expansion begins $[2n,m,\ldots]$, and off these intervals $f_2(\theta)=1$:
\begin{align*}
\int_{S^1} \log f_2(\theta)d \theta & = \sum_{n,m=1}^{\infty} \int_{\frac{m}{2nm+1}}^{\frac{m+1}{2n(m+1)+1}}\log m d\theta\\
&= \sum_{n,m=1}^{\infty} \frac{\log m}{(2nm+1)(2n(m+1)+1)}\\
&< \infty.
\end{align*}
Finally, note that if $\{f_2(\theta_i)\}$ is exceptionally irregular, then the ergodic average of $\log(f_2)$ cannot converge to a positive real number, while as the integrable function $\log(f_2)$ is positive on a set of positive measure, the ergodic theorem guarantees such convergence almost surely.
\end{proof}
\end{lemma}
\begin{lemma}
$\log(f_1) \in L^1(S^1,d\mu_g)$.
\begin{proof}
The proof is similar to the computation in Lemma \ref{lemma - no fast growth}, using both \eqref{eqn - estimate even} and \eqref{eqn - estimate odd} to show that
\[\int_{0}^1 \left| \log(f_1 (x)) \right| dx < \sum_{k=1}^{\infty} \frac{\log(2k+1)}{(2k+1)(2k+2)} + \sum_{n,m=1}^{\infty} \frac{\log(2n+1)+\log(m+1)}{(2nm+1)(2n(m+1)+1)}. \qedhere\]
\end{proof}
\end{lemma}

Combining Lemma \ref{lemma - no fast growth} and Proposition \ref{proposition - hausdorff is box}, for almost every $\theta$ we have
\[\dim_H{H_{\theta}} = \liminf_{n \rightarrow \infty} \frac{\sum_{i=0}^{n-1} \log(f_2(\theta))}{-\sum_{i=0}^{n-1} \log(f_1(\theta))}.\]  By the ergodic theorem, however, for almost every $\theta$ the right side converges as a proper limit to the ratio
\[c=\frac{\int_{S^1} \log(f_2(\theta))d\mu_g}{-\int_{S^1} \log(f_1(\theta))d\mu_g}.\]
That $c \neq 0$ is seen by noting again that the nonnegative function $\log(f_2)$ is in fact positive on a set of positive measure.  By verifying the pointwise inequality $-\log(f_1)>\log(f_2)$ for $\theta<1/2$ (using \eqref{eqn - estimate even} for $a_1$ even and $\log(f_2)=0$ for $a_1$ odd), we see $c \neq 1$ (both $f_1(\theta)$ and $f_2(\theta)$ equal one for $\theta>1/2$).

Note that we have proved something slightly stronger than the original statement of Theorem \ref{theorem - heavy constant}.  The set $H_{\theta}$ is always comprised of a (possibly empty) set of isolated points $S$ together with a set $E$, where $E$ is almost surely of both Hausdorff and box dimension $c$, where 
\[0<c=\frac{\int_{S^1} \log(f_2(\theta))d\mu_g}{-\int_{S^1} \log(f_1(\theta))d\mu_g}<1.\]
We also see that $E$ is a perfect set (closed with no isolated points) for all $\theta$ such that $f_2(\theta_i)\neq 1$ infinitely many times, which happens almost surely (via ergodicity of $g$ and continuity of $\mu_g$).  Note that $H_{\theta}$ is \textit{exactly} this perfect set exactly when $S = \emptyset$, which happens exactly when $a_1(\theta_i)=0\bmod 2$ for all $i$.  These are exactly the $\theta$ for which all partial quotients of odd index are themselves even, which we have remarked are exactly the $\theta$ for which $H^*_{\theta}=\{0\}$:

\begin{corollary}
$H_{\theta}$ is a perfect set if and only if $H^*_{\theta}=\{0\}$, if and only if $a_{2i+1}(\theta)=0\bmod 2$ for all $i=0,1,\ldots$
\end{corollary}

\section{Proof of Theorem \ref{theorem - heavy dense}}
Let $d \in [0,1]$; we begin by showing that there is some $\theta$ such that $\dim_H(H_{\theta})=d$.  As the general proof is simplified by the assumption $d \neq 0$, let us take care of that particular case by hand.  Let
\[\theta = [2!, 1, 3!, 1, 4!, 1, \ldots],\]
so $\theta_i = [(i+2)!,1,(i+3)!,1,\ldots]$.  As the sequence $f_2(\theta_i)$ is not exceptionally irregular (it is constantly two), we have (using \eqref{eqn - estimate even})
\begin{align*}
\textrm{dim}_H(H_{\theta}) & = \liminf_{n \rightarrow \infty} \frac{n \log 2}{-\sum_{i=0}^{n-1} \log (1-(i+2)!\theta_i)}\\
&\leq \liminf_{n \rightarrow \infty} \frac{n \log 2}{\sum_{i=0}^{n-1} \log ((i+2)!)},
\end{align*}
which is clearly zero.  Note that $H_{\theta}$ is a perfect set and therefore uncountable.

So now let us proceed on the assumption that $d \neq 0$.  Define $m_i = 2^i$ and 
\[n_i = 1+ \integer{\frac{(m_i + 1)^{\frac{1}{d}-1}}{2}}.\]
By construction,
\begin{equation}\label{equation - goes to d}
\frac{\log(m_i+1)}{\log(2n_i)+\log(m_i+1)} \rightarrow d.
\end{equation}
That the values $\log(m_i/(m_i+1))$ vanish exponentially in $i$, and $\log(2n_i/(2n_i+1))$ are either constant (if $d=1$) or vanish exponentially in $i$ (otherwise) further guarantees that
\[\frac{1}{\sum_{i=0}^n \log (m_i+1)} \left( \sum_{i=0}^n \log \left(\frac{2n_i}{2n_i+1}\right) + \sum_{i=0}^n \log \left( \frac{m_i}{m_i+1}\right) \right) \rightarrow 0,\]
which is seen to imply
\begin{equation}\label{equation - n chosen correctly}
\liminf_{n \rightarrow \infty} \frac{\sum_{i=0}^{n}\log(m_i+1)}{\sum_{i=0}^n \log(2n_i m_i)} = \liminf_{n \rightarrow \infty} \frac{\sum_{i=0}^n \log (m_i+1)}{\sum_{i=0}^n \log\left((2n_i+1)(m_i+1)\right)}.
\end{equation}

Let $\theta = [2n_0,m_0,2n_1,m_1,\ldots]$, so $\theta_i=[2n_i,m_i,2n_{i+1},m_{i+1},\ldots]$.  As the $m_i$ are not exceptionally irregular \eqref{eqn - not exceptionally irregular}, we have
\[\textrm{dim}_H(H_{\theta}) = \liminf_{n \rightarrow \infty} \frac{\sum_{i=0}^n \log(m_i+1)}{-\sum_{i=0}^n \log(1-2n_i \theta_i)}.\]
The denominator may be approximated from above by $\log\left((2n_i+1)(m_i+1)\right)$ and from below by $\log(2n_i m_i)$ using \eqref{eqn - estimate even}:

\[ \liminf_{n \rightarrow \infty} \frac{\sum_{i=0}^n \log (m_i+1)}{\sum_{i=0}^n \log\left((2n_i+1)(m_i+1)\right)}  \leq \textrm{dim}_H(H_{\theta}) \leq \liminf_{n \rightarrow \infty} \frac{\sum_{i=0}^n \log(m_i+1)}{\sum_{i=0}^n \log(2n_i m_i)},\]
and by applying \eqref{equation - n chosen correctly}, the two sides have the same limit infimum, and by \eqref{equation - goes to d} this common value is $d$.

As $\dim_H(H_{\theta})$ is a $g$-invariant function of $\theta$ for \textit{every} irrational $\theta$, the theorem now follows from the fact that
\[\bigcup_{i=0}^{\infty} g^{-i}(\theta)\] is a dense set for every irrational $\theta$; given any finite string of initial partial quotients $a_i(\theta)$, a sufficiently high power of $g$ will have truncated all of them from the continued fraction expansion.

\section{Examples}

\begin{example}\label{example - best heavy}
For $\theta = \sqrt{2} \bmod 1=[2,2,2,\ldots]$, $H_{\theta}$ is a perfect set of both Hausdorff and box dimension
\[\frac{\log(3)}{\log(3+2\sqrt{2})} \approx .623,\] and $H^*_{\theta}=\{0\}$.
\begin{proof}
Since $\theta_i = [2,2,2,\ldots]$ for all $i$, we always form $E_{i+1}$ from $E_i$ by considering $a_2+1=3$ subintervals scaled in length by $1-2\theta = 3-2\sqrt{2}=(3+2\sqrt{2})^{-1}$.  So the Hausdorff and box dimensions of this simple Cantor set are both given by
\[ \dim_H(H_{\theta})=\dim_B(H_{\theta})  = \frac{\log 3}{\log(3+2 \sqrt{2})}.\]
We have already remarked that if all $a_{2i+1}(\theta)=0\bmod 2$ then $H^*_{\theta}=\{0\}$.
\end{proof}
\end{example}

It is tempting to think that the tail of the continued fraction expanion of $\theta$ is what controls the development of the heavy set.  This is not the case, however: it is dependent on orbit of $\theta$ under $g$, which might be significantly different than the orbit under the Gauss map $\gamma$.

\begin{example}\label{example - worst heavy}
Let $\theta$: $\sqrt{2}/{2}=[1,2,2,2,\ldots]$.  Then $H_{\theta}$ consists of countably many isolated points and exactly one accumulation point, given by $H^*_{\theta}=\{\theta/2\}$.
\begin{proof}
In this case, $a_1(\theta_i)\in \{1,3\}$ for all $i$, so $E_i=E^*_i$ for all $i$.  So the set $E$ is in fact the singleton $H^*_{\theta}$; every other point in $H_{\theta}$ (which must be infinite by Corollary \ref{corollary - always infinite}) therefore belongs to the set of isolated points $S$.  Note that even though $a_i(\theta)=2$ for all $i \neq 1$, $a_1(\theta_i)\neq 2$ for all $i$.

The accumulation point of $H_{\theta}$ is given by $\cap E_i$, but for this choice of $\theta$ we have $E_i=E^*_i$, so the accumulation point is $H^*_{\theta}$.  To compute the value of this point, note that $a_1(\theta_i)\neq 1$ if and only if 
\[\theta_i=[3,2,2,\ldots]=1/(2+\sqrt{2}),\] for which $\delta_i = 1-2\theta_i=\sqrt{2}-1$.  Using Table \ref{table - decomposition} to construct the sets $E_i=E^*_i$, we see that as every other $a_1(\theta_i)=1$, we alternate placing $E_{i+1}$ at the top or bottom of $E_i$, each time scaled by $\delta=\sqrt{2}-1$:
\[E_0=E_1=[0,1/2], \quad E_2=E_3=[1/2-\delta/2,1/2], \quad E_4=E_5=[1/2-\delta,1/2+\delta^2/2], \ldots,\]
and the intersection is given by the geometric series
\[H^*_{\theta} = \frac{1}{2}\sum_{i=0}^{\infty}(-1)^i(\sqrt{2}-1)^i = \frac{1}{2(1+(\sqrt{2}-1))}=\frac{1}{2\sqrt{2}}=\frac{\theta}{2}.\qedhere \]
\end{proof}
\end{example}

\begin{example}\label{example - e and related heavy}
As a final pair of examples of general interest, consider 
\[\begin{array}{r c l}
\alpha 	&=e-2 	&= [1,2,1,1,4,1,1,6,1,\ldots,1,2k,1,\ldots] \\
\beta 	&=\frac{e-1}{e+1} &= [2, 6, 10, 14, 18, \ldots, 4k-2,\ldots].
\end{array}\]
Then $\dim_H(H_{\alpha})=0$, $\dim_H(H_{\beta}) = 1/2$.

The value $\beta$ is easier to consider: $\beta_i=[8i+2, 8i+6, \ldots ]$.  The sequence $f_2(\beta_i)$ is not grow exceptionally irregular (it is an arithmetic sequence), so
\[\dim_H(H_{\beta}) = \liminf_{n \rightarrow \infty} \frac{\sum_{i=0}^{n-1} \log(8i+7)}{-\sum_{i=0}^{n-1}\log(1-(8i+2)\theta_i)},\] which is bounded above and below (using \eqref{eqn - estimate even}) by
\[\begin{split}\liminf_{n \rightarrow \infty} \frac{\sum_{i=0}^{n-1} \log (8i+7)}{\sum_{i=0}^{n-1}\log(8i+2) + \sum_{i=0}^{n-1}\log(8i+6)}\geq \dim_H(H_{\beta})\\ \geq \liminf_{n \rightarrow \infty} \frac{\sum_{i=0}^{n-1} \log (8i+7)}{\sum_{i=0}^{n-1}\log(8i+3) + \sum_{i=0}^{n-1}\log(8i+7)}.\end{split}\]
Both sides are seen to have proper limit one half.  As all $a_1(g^i \beta)=0\bmod 2$, there are no isolated points in $H_{\beta}$.

Now considering $\alpha$, we must examine how $H_{\alpha}$ decomposes into the sets $E_i$ based on $\alpha_i$.  
\begin{table}
\begin{tabular}{| c | c |  c | c|}
\hline $i$ & $\alpha_i$ & $\log(f_2(\alpha_i))$ & $-\log(f_1(\alpha_i))\geq $ \\
\hline  \hline $5k$ & $[1,4k+2,1,\ldots]$ &$0$ &$0$\\
\hline $5k+1$ & $[4k+3,1,1,4k+4,\ldots]$ & $0$&$\log(2k+2)$\\
\hline $5k+2$ & $[1,1,1,4k+4,\ldots]$ & $0$&$0$\\
\hline $5k+3$ & $[2,1,4k+4,\ldots]$ & $\log 2$&$\log 2$\\
\hline $5k+4$ & $[4k+4,1,1,4k+6,\ldots]$ & $\log 2$&$\log(4k+4)$\\
\hline 
\end{tabular}
\caption{The $5$-step pattern in construction $H_{e}$.  Bounds for $f_1$ are established using \eqref{eqn - estimate even}, \eqref{eqn - estimate odd}.}\label{table - build He}
\end{table}
Referring to Table \ref{table - build He}, we see how the sequence $\alpha_i$ develops with a $5$-step pattern with respect to the functions $f_1$ and $f_2$.  We directly compute the lower box dimension of the set of non-isolated points in $H_{\alpha}$:
\[\dim_H{H_{\alpha}}\leq \dim_{b}(E) = \liminf_{n \rightarrow \infty} \frac{\sum_{k=0}^{n-1} (2\log 2)}{\sum_{k=0}^{n-1} \left( \log(2k+2) + \log(2) + \log (4k+4) \right)}=0.\]
\end{example}

\section*{Acknowledgements}
The author is supported by the Center for Advanced Studies at Ben Gurion University of the Negev as well as the Israeli Council for Higher Education.
\bibliography{../masterbib}{}

\begin{thebibliography}{1}

\bibitem{MR2515388}
Michael Boshernitzan and David Ralston.
\newblock Continued fractions and heavy sequences.
\newblock {\em Proc. Amer. Math. Soc.}, 137(10):3177--3185, 2009.

\bibitem{MR2118797}
Kenneth Falconer.
\newblock {\em Fractal geometry}.
\newblock John Wiley \& Sons Inc., Hoboken, NJ, second edition, 2003.
\newblock Mathematical foundations and applications.

\bibitem{springerlink:10.1007/BF01896805}
G.~Hal\'{a}sz.
\newblock Remarks on the remainder in {B}irkhoff's ergodic theorem.
\newblock {\em Acta Mathematica Hungarica}, 28:389--395, 1976.
\newblock 10.1007/BF01896805.

\bibitem{MR931186}
Yuval Peres.
\newblock A combinatorial application of the maximal ergodic theorem.
\newblock {\em Bull. London Math. Soc.}, 20(3):248--252, 1988.

\bibitem{ralston}
David Ralston.
\newblock Substitutions and $1/2$-discrepancy sums of $\{x+n \theta\}$.
\newblock http://arxiv.org/abs/1105.5810v1.

\bibitem{ralston2}
David Ralston.
\newblock Substitutions and $1/2$-discrepancy sums of $\{x+n \theta\}$ \rm{II}.
\newblock http://arxiv.org/abs/1105.6301v1.

\end{thebibliography}
\bibliographystyle{plain}
\end{document}